\documentstyle[12pt]{article}
\title{Supersimple fields and division rings}
\author{A. Pillay\thanks{Partially supported by NSF grant DMS
96-96268}\\University of Illinois and MSRI \and T.
Scanlon\thanks{Supported by an NSF Postdoctoral Fellowship}\\MSRI \and F.O.
Wagner\thanks{Supported by DFG grant Wa899/2-1}\\University of Oxford and MSRI}
\newtheorem{Theorem}{Theorem}[section]
\newtheorem{Proposition}[Theorem]{Proposition}
\newtheorem{Lemma}[Theorem]{Lemma}
\newtheorem{Fact}[Theorem]{Fact}

\newtheorem{Remark}[Theorem]{Remark}
\newtheorem{Corollary}[Theorem]{Corollary}

\begin{document}
\maketitle
\begin{abstract}
It is proved that any supersimple field has trivial Brauer group, and more
generally that any supersimple division ring is commutative. As prerequisites
we prove several results about generic types in groups and fields whose theory
is simple.
\end{abstract}
\section{Introduction}
Simple theories were introduced by Shelah in \cite{Shelah}. In \cite{Kim}, Kim,
continuing Shelah's work, showed how the theory of forking transfers almost
completely from the stable context to the more general simple context. In
\cite{Kim-Pillay} the ``Independence Theorem" was proved for simple theories,
giving a satisfactory analogue of the theory of stationarity from stability
theory.

In \cite{Pillay}, Poizat's theory (\cite{Poizat}) of generic types and
stabilizers in stable groups was generalized to the case of groups definable in
simple theories. Further generalizations appear in \cite{Wagner}.

Stability-theoretic algebra studied, among other things, the algebraic
consequences of imposing stability-theoretic conditions on a group, ring or
field. Among the main results in the area was the Macintyre-Cherlin-Shelah
theorem (\cite{Macintyre}, \cite{Cherlin-Shelah}) saying that any infinite field
whose theory is superstable must be algebraically closed.

Supersimple fields (namely fields with supersimple theory) on the other hand
form a broader class. In \cite{Hrushovski} it is shown that any perfect
pseudo-algebraically closed ($PAC$) field with ``small" absolute Galois group is
supersimple; in fact of $SU$-rank $1$ in the language of rings. (Ultraproducts
of finite fields furnish examples.) The work in this paper is partly
motivated by
the conjecture that these are the only cases. It follows from
\cite{Pillay-Poizat} that any supersimple field is perfect and has small
absolute Galois group. In the present paper we show that a supersimple field $K$
has trivial Brauer group. As supersimplicity is preserved under finite
extensions, this is, by
\cite{Serre}, X.Prop.11, equivalent to the norm map $N : L^{*} \rightarrow
K^{*}$ being surjective for any finite Galois extension $L$ of $K$. Triviality
of the Brauer group also implies that any rational variety over $K$ has a
$K$-rational point, yielding in a sense a first approximation to the conjecture
that any supersimple field must be $PAC$.

We assume acquaintance with the notions and machinery from \cite{Kim-Pillay} and
\cite{Pillay}. In section 2, some additional results are obtained about generic
types and connected components of groups in simple theories. In section 3, we
show that if $K$ is a field in a simple theory then the notions of additive and
multiplicative generic coincide. We also make some observations about the
interaction between the additive and multiplicative connected components. These
results are applied in section 4 to show triviality of the Brauer group in the
supersimple case. Some additional algebraic arguments yield in section 5
commutativity of any supersimple division ring.

The first author would like to thank Zoe Chatzidakis for many fruitful
discussions around the time when he first began considering the issues
dealt with
in this paper. In particular the result that if $K$ is a supersimple field of
characteristic $0$, then any rational {\em curve} over $K$ has a $K$-rational
point, was obtained together with her.

The result on the triviality of the Brauer group was obtained by all three
authors, and the results in section 5 on division rings were obtained by Wagner.

\section{More on generic types and stabilizers}
Let us fix a saturated model $\bar{M}$ of a simple theory, and an infinite
group $G$ which is type-definable in $\bar{M}$ over $\emptyset$ say. All
complete types $p(x)$ we consider will be types of elements of $G$ (namely we
assume $``p(x) \rightarrow x \in G"$). We will often wish to work with Lascar
strong types (so that we can use the Independence Theorem). This is most easily
accomplished by working with types over models. So $M$ will denote a small
elementary substructure of $\bar{M}$, which may vary. $A,B$ as usual denote
small subsets of $\bar{M}$.

We briefly recall (from \cite{Pillay}) the notions of generic type,
stabilizer, and connected component,

\noindent
{\bf Generics.} Let $a\in G$. Then $tp(a/A)$ is (left-) generic (for $G$) if
whenever $b\in G$ is independent from $a$ over $A$ then $b.a$ is independent
from $A\cup \{b\}$ over $\emptyset$. (We also say that $a$ is a generic element
of $G$ over $A$.) We have associated notions of a (type)-definable subset of $G$
being generic in $G$. Moreover the notions of left-generic and right-generic
coincide.

\noindent
{\bf Stabilizers.}  For $p(x)\in S(M)$, $St(p)$ is defined to be $\{b\in
G$: for some realization $a$ of $p$ which is independent from $b$ over $M$,
$b.a$ realises $p$ (and is independent from $b$ over $M$\}.$St(p)$ is
type-definable over $M$. Let $Stab(p)$ be the subgroup of $G$ generated by
$St(p)$. Then it turns out that $Stab(p)$ is type-definable over $M$, and
moreover that $St(p)$ is ``large" in $Stab(p)$ in the sense that any generic
element of $Stab(p)$ over $M$ is in $St(p)$.

\begin{Remark} One can also define $St(p,q)$ for $p(x), q(x) \in S(M)$ in the
obvious way: $St(p,q)$ = $\{b\in G$: for some $a$ realising $p$ such that $a$
is independent from $b$ over $M$, $b.a$ realises $q$ and is independent from
$b$ over $M$\}. $St(p,q)$ is again type-definable over $M$ and one sees easily,
using the Independence Theorem, that $St(p,q)$ is, if nonempty a large subset of
a translate of $Stab(p)$.
\end{Remark}
{\bf Connected components.} For any set $A$ of parameters, $G^{0}_{A}$, the
connected component of $G$ over $A$, is by definition the smallest
type-definable over $A$ subgroup of $G$ of bounded index. A key fact is that
$p(x)\in S(M)$ is a generic type of $G$ if and only if $Stab(p)$ contains
$G^{0}_{M}$. Also clearly any generic type of $G^{0}_{A}$ is a generic type of
$G$. In contrast to the stable situation, the connected component of $G$
over $A$ may vary with $A$. On the other hand $G^{0}_{A}$ is always a normal
subgroup of $G$.

If $G$ happens to be stable then we know that $G^{0}$ has a
unique generic type, whose stabilizer is $G^{0}$. The next lemma gives an
analogue of this in the simple case.
\begin{Proposition} Let $p(x), q(x), r(x) \in S(M)$ be generic types of
$G^{0}_{M}$. Then there are realizations $a,b,c$ of $p,q,r$ respectively, which
are pairwise $M$-independent, and with $a.b = c$.
\end{Proposition}
{\em Proof.} This can be seen using Remark 2.1, but we give a direct proof. Let
$b, c$ be independent (over $M$) realizations of $q,r$ respectively. Let $a'$
= $c.b^{-1}$. Then $a'.b = c$, $a'$ is generic in $G^{0}_{M}$ over $M$ and
$\{a',b,c\}$ is pairwise $M$-independent. We can similarly find $d$ generic
in $G^{0}_{M}$ over $M$, such that $a = a'.d$ realises $p$ and $\{d,a',a\}$
is pairwise $M$-independent. By the facts above on connected components and
generics, $Stab(q)$ contains $G^{0}_{M}$ and hence any generic element of
$G^{0}_{M}$ over $M$ is in $St(q)$. In particular $d^{-1} \in St(q)$.  By the
Independence Theorem we can assume that $\{d,a',b\}$ is pairwise
$M$-independent and that $d^{-1}.b = b'$ realises $q$. Then $(a'.d).b' = c'$
realises $r$. That is $a.b' = c$, and easily $\{a,b',c\}$ is pairwise
$M$-independent.

\section{Simple fields}
In this section $F$ will
be an infinite field which is type-definable over $\emptyset$ in the saturated
model $\bar{M}$ of a simple theory $T$. Types $p(x)$, $q(x)$ will be types of
elements of $F$.
We are going to apply Proposition 2.2 in order to understand the interaction
between generic additive subgroups and generic multiplicative subgroups in
simple
fields. (In the stable case the situation is clear: a stable field is connected
both additively and multiplicatively.) We begin by pointing out that in simple
fields additive and multiplicative generics coincide.

\begin{Proposition} Let $F$ be an infinite field which is type-definable in
$\bar{M}$ over $\emptyset$. Then a type $p(x)\in S(A)$ is an additive generic
type of $F$ iff it is a multiplicative generic type of $F$.
\end{Proposition}
{\em Proof.} We begin with
\newline
\noindent
{\em Claim 1.} If $p(x) \in S(A)$ is additively generic then for any
nonzero $a\in F$ and realisation $c$ of $p$ independent from $a$ over $A$,
$tp(a.c/A\cup\{a\})$ is an additive generic of $F$.
\newline
\noindent
{\em Proof of Claim 1.} Let $a\in F$ and $M$ contain $A\cup\{a\}$. Let $p'$
be a nonforking extension of $p$ over $M$, realised by $b$. Then $p'$ remains
an additive generic of $F$. So the additive stabilizer of $p'$, $Stab^{+}(p')$
is a type-definable over $M$ subgroup of $F^{+}$ of bounded index. It follows
that $Stab^{+}(a.p')$ = $a.Stab^{+}(p')$ also has bounded index. By the remarks
above on connected components $a.p'$ is an additive generic type of $F$. In
particular $a.p'$ does not fork over $\emptyset$.
\newline
\noindent
From Claim 1, using the fact that $tp(a.c/A\cup\{a\})$ there, being an additive
generic, does not fork over $\emptyset$, we see that:
\newline
\noindent
{\em Claim 2.} Any additive generic type is a multiplicative generic type.
\newline
\noindent
{\em Claim 3.} Any multiplicative generic type of $F$ is an additive generic
type of $F$.
\newline
\noindent
{\em Proof of Claim 3.} Let $q(x)\in S(A)$ be a multiplicative generic type. We
may replace $q$ any time by a nonforking extension. By Claim 2 (and the
existence of generic types) there is an additive generic type $p(x) \in S(A)$
which is also multiplicatively generic. It follows that there is $c \in F$, $M$
containing $A\cup\{c\}$ and nonforking extensions $q',p'$ of $q,p$ over $M$ such
that $c.q'$ = $p'$. So $c.q'$ is an additive generic type. By Claim 1, $q'$ is
also an additive generic type.

\begin{Lemma} Let $T$ be a multiplicative subgroup of $F$
of bounded index, type-definable over $A$. Then every nonzero coset $a.T$ of $T$
meets $(F^{+})^{0}_{A}$ in a generic set.
\end{Lemma}
{\em Proof.} Note that if  $B\supset A$, then $(F^{+})^{0}_{B}$ is contained in
$(F^{+})^{0}_{A}$. So we may assume that $A$ is a model $M$ which contains a
representative of each coset $c.T$ of $T$ in $F^{*}$. Note that
\vspace{2mm}
\noindent
(*) For each nonzero $c \in F^{M}$,
$c.(F^{+})^{0}_{M}$ = $(F^{+})^{0}_{M}$ (as multiplication by $c$ is an additive
endomorphism of $F$ which interchanges $M$-definable sets).
\vspace{2mm}
\noindent
Let $p(x)$ be a
generic type of $G^{0}_{M}$ over $M$. Let $S$ be a coset of $T$ in $F^{*}$.
Then there is $c\in F^{M}$ such that $c.p$ is the type of an element of
$S$. By Claim 1 in the proof of Proposition 3.1, together with (*), $c.p$ is a
generic type of $(F^{+})^{0}_{M}$. This proves the lemma.

\begin{Remark} In fact, with the assumptions of Lemma 3.2, every coset of
$(F^{+})^{0}_{A}$ in $F^{+}$ meets every coset of $T$ in $F^{*}$ in a generic
set.
\end{Remark}
{\em Proof.} We may assume that $A$ is a model $M$ which contains
representatives of all cosets of $T$ in $F^{*}$. To prove the
remark, it is enough (using Lemma 3.2) to show that if $S$ is a coset of
$(F^{+})^{0}_{M}$ in $F^{+}$ then $S\cap T$ is generic. Let $b\in S$ be nonzero.
Let $M'$ be a model containing $M \cup\{b\}$. By Lemma 3.2, we can find generic
$d$ in $(F^{+})^{0}_{M'}$ over $M'$ such that $d \in b^{-1}.T$. Let $c = d-1$.
Then $c$ is generic in $(F^{+})^{0}_{M'}$ over $M'$ too, and as in the proof of
Lemma 3.2, $b.c$ is generic in $(F^{+})^{0}_{M'}$ over $M'$. So $b.d$ =
$b.c + b$
is generic in $S$ over $M$. On the other hand $b.d \in T$. Thus $S \cap T$ is
generic, as required.

\begin{Proposition} Let $T$ be any multiplicative subgroup of $F$ of bounded
index. Let $S_{1}, S_{2}$ be cosets of $T$ in $F^{*}$. Then $S_{1} + S_{2}$
(the set of $s_{1} + s_{2}$ for $s_{1} \in S_{1}$ and $s_{2} \in S_{2}$)
contains $F^{*}$. In fact for any nonzero $d\in F$ there are $a,b \in S_{1},
S_{2}$ respectively, each generic over $d$ such that $a+b = d$.
\end{Proposition}
{\em Proof.} Again we work over a reasonably saturated model $M$ (so in
particular $T, S_{1}, S_{2}$ are defined over $M$). Let $d\in F^{M}$ be
nonzero, and let $S_{3}$ be $d.T$. By Lemma 3.2, let $p(x)$, $q(x)$, $r(x)$ be
generic types of $(F^{+})^{0}_{M}$ such that $p(x) \rightarrow x\in S_{1}$,
$q(x)
\rightarrow x\in S_{2}$ and $r(x) \rightarrow x\in S_{3}$. By Proposition 2.2,
there are realisations $a,b,c$ of $p,q,r$ respectively, such that $a+b = c$,
where moreover $\{a,b,c\}$ is pairwise independent over $M$. Multiplying on the
right by $c^{-1}.d \in T$, we obtain $a' \in S_{1}$, $b'\in S_{2}$ with $a' + b'
= d$. Clearly each of $a', b'$ are generic over $M$. This completes the proof.

\vspace{5mm}
\noindent
[Here is a sketch of another proof of 3.4, avoiding use of Proposition 2.2
and containing some possibly useful ideas: Work again over a model $M$. Let
$X$ =
$St(S_{1},S_{2})$ which by definition equals $\{a\in F:$ for each generic type
$p(x)\in S(M)$ of $S_{1}$ there is $c$ realising $p$ independent from $a$ such
that $a.c \in S_{2}$\}. $X$ is type-definable over $M$ and moreover (see 2.2)
there is a type-definable set $X'$, a union of cosets of $(F^{+})^{0}_{M}$ in
$F^{+}$ such that $X \subset X'$ and also every generic element of $X'$ over $M$
is in $X$. Note that $X$ is invariant under multiplication by $T$. Let now $H$
be the set-theoretic stabilizer of $X'$. $H$ is a type-definable over $M$
subgroup of $F^{+}$ containing $(F^{+})^{0}_{M}$. $H$ is also invariant under
multiplication by $T$, but we know that $(F^{+})^{0}_{M}$ intersects each coset
of $T$ in $F^{*}$. Thus $H$ = $F^{+}$. It follows that $X'$ = $F$, and thus $X$
contains all generic types over $M$. As $X$ is also invariant under $T$ it
follows that $X$ contains $F^{*}$.]

\vspace{5mm}
\noindent
The characteristic $0$ case of the following was observed earlier together with
Zoe Chatzidakis.
\begin{Corollary} If $F$ is supersimple, then any conic defined over $F$ has an
$F$-rational point.
\end{Corollary}
{\em Proof.} A conic defined over $F$ can be put into the form $x^{2}+ay^{2}
= b$ for some nonzero $a,b \in F$. As we note in the next section,
supersimplicity of $F$ implies that the squares form a definable subgroup of
$F^{*}$ of finite index. By 3.4, the equation has a solution in $F$.

\section{Triviality of the Brauer group for
supersimple fields}
We will prove that if $F$ is an infinite field possibly with
extra structure, whose theory is supersimple, then $F$ has trivial Brauer group.
We will see in the course of the proof that this is a first order property, so
we will assume $F$ to be saturated. We will work with the formalism of the
previous section, namely we assume $F$ to be definable over $\emptyset$ in a big
model $\bar{M}$ of a supersimple theory.

At this point it is worth mention the relevant facts about the $SU$-rank,
most of which will be used in the next section. Each complete type has
ordinal valued $SU$-rank, and the generic types of a
(type-)definable group $G$ are precisely the types of maximal $SU$-rank in $G$
(which exist). The following is mentioned explicitly in \cite{Wagner}

\begin{Fact} (i) for any elements $a,b$ and set $A$ of parameters,
$SU(a/A\cup\{b\}) + SU(b/A) \leq SU(a,b/A) \leq SU(a/A\cup \{b\}) \oplus
SU(b/A)$.
\newline
\noindent
(ii) Let $G$ be a type-definable group and $H$ a type-definable subgroup, then
$SU(H) + SU(G/H) \leq SU(G) \leq SU(H) \oplus SU(G/H)$
\newline
\noindent
(iii) Suppose $D$ is a type-definable division ring, then $D$ is definable and
has monomial $SU$-rank.
\end{Fact}

The relevant consequences for this section are:
\begin{Fact} (i)
for each n, the group of nth. powers $(F^{*})^{n}$ has finite index in $F^{*}$,
\newline
\noindent
(ii) $F$ is perfect.
\end{Fact}
{\em Proof.} (i) Let $a \in F$ be generic over $\emptyset$. Let $b =
a^{n}$. Then
$a \in acl(b)$, so $SU(a/b) = 0$. Thus by Fact 4.1. (i) $SU(a) = SU(b)$, whereby
$b$ is generic in $F$ too. Thus $(F^{*})^{n}$ is a generic subgroup of $F^{*}$
and  has finite index (as it is definable).
\newline
\noindent
(ii) Supposing $F$ to have characteristic $p > 0$, note that as in (i) the
additive group of $F^{p}$ has finite index in $F$. But $F/F^{p}$ is a vector
space over the field $F^{p}$, a contradiction (the latter being infinite)
unless $F^{p} = F$.

\noindent
Note that any finite extension of $F$ is also definable in $\bar{M}$, so Fact
4.2 applies to all finite extensions of $F$.

We will now recall relevant notions and facts concerning the Brauer group,
the norm map and Galois cohomology. In fact one can extract a general result
which says that for the Brauer group of every finite extension of a field $F$ to
be trivial it is enough that for any finite extension K of F and Kummer
extension $L$ of $K$, $N_{L/K}$ is surjective. (Here it is assumed that every
finite extension of $F$ is perfect.) Proposition 3.4 will allow us to conclude
triviality of the Brauer group for supersimple $F$. Rather than simply
state this reduction, we will include an explanation of it as part of the
proof, which will entail giving some definitions.

We first discuss the Brauer group. The reader can look at \cite{Jacobson} and
\cite{Serre} for further details. We will assume that $F$ and all its finite
extensions are perfect.

By a central simple algebra over $F$ we mean a finite
dimensional $F$-algebra $A$ whose centre is $F$ and which has no nontrivial
two-sided ideals. If $A$ and $B$ are two such objects then so is the tensor
product $A \otimes _{F} B$. $A$ and $B$ are called equivalent (or
similar) if for some $m,n$ the matrix algebras $M_{m}(A)$, $M_{n}(B)$ are
isomorphic (as $F$-algebras). The tensor product operation respects this
equivalence relation and turns the set of classes into an abelian group. This
group is called the {\em Brauer group of $F$}, $Br(F)$, and is an important
invariant of the field $F$. Any central simple $F$-algebra $A$ will be a matrix
algebra over a certain finite dimensional division algebra $D$ with centre $F$.
Moreover the equivalence class of $A$ is determined by and determines the
isomorphism type of $D$. The trivial element of $Br(F)$ then corresponds to $F$
itself. On the other hand, for any central simple algebra $A$ over $F$ there
i       s some finite extension $K$ of $F$ such that $A \otimes K$ is isomorphic
to a matrix algebra over $K$, hence represents the trivial element of $Br(K)$.
($K$ is called a splitting field for $A$.)

If $K$ is a finite extension of $F$ then tensoring with $K$ determines a
homomorphism from $Br(F)$ into $Br(K)$. The kernel of this homomorphism is
denoted $Br(K/F)$. From the previous paragraph, $Br(F)$ is the union of
all $Br(K/F)$ as $K$ runs over finite extensions of $F$, in fact, over finite
Galois extensions of $F$.

For $K$ a finite Galois extension of $F$ with Galois group $G$, there is a
classical isomorphism of $Br(K/F)$ with the Galois cohomology group
$H^{2}(G,K^{*})$: given a $2$-cocycle $f:G \times G \rightarrow K^{*}$, define
a $F$-algebra structure on the $K$-vector space $A$ with basis $\{u_{s}:s\in
G\}$ by: $u_{s}\cdot u_{t} = f(s,t)u_{st}$ and for $e\in K$, $u_{s}e =
(se)u_{s}$. $A$ then becomes a central simple $F$-algebra $A_{f}$, and the map
$f \rightarrow A_{f}$ determines an isomorphism between $H^{2}(G,K^{*})$ and
$Br(K/F)$.
\vspace{3mm}
\noindent
We have the following fact (see \cite{Serre}, X.6):
\begin{Fact} Let $F <K <L$ where both $K$ and $L$ are Galois extensions of $F$
with Galois groups $H, G$ respectively. Then there is an exact sequence
\newline
\noindent
$0\rightarrow H^{2}(H,K^{*}) \rightarrow H^{2}(G,L^{*})\rightarrow
H^{2}(Gal(L/K),L^{*})$.
\end{Fact}
\vspace{3mm}
\noindent
We also need (\cite{Serre}, IX):
\begin{Fact} Let $G$ be a finite group and $A$ a $G$-module. Let $n\geq 1$.
Suppose that for all primes $p$, $H^{n}(G_{p},A) = 0$ where $G_{p}$ is a Sylow
subgroup of $G$. Then $H^{n}(G,A) = 0$.
\end{Fact}

Finally, for $K$ a finite Galois extension of $F$, the norm map
$N_{K/F}:K^{*}\rightarrow F^{*}$ is the map which takes any $a\in K^{*}$ to the
product of all $sa$ where $s$ runs over $Gal(K/F)$. See Theorem 8.14 of
\cite{Jacobson} for:
\begin{Fact} Let $K$ be a cyclic extension of $F$ with
Galois group $G$. Then $H^{2}(G,K^{*})$ is isomorphic to the quotient group
$F^{*}/N_{K/F}(K^{*})$.
\end{Fact}

We can now prove:
\begin{Theorem} ($F$ a supersimple field.) $Br(F)$ is trivial.
\end{Theorem}
{\em Proof.} We will prove by induction on $n$, that for every finite extension
$K$ of $F$ and every Galois extension $L$ of $K$ of degree $n$, $Br(L/K)$ (or
equivalently $H^{2}(Gal(L/K),L^{*})$) is trivial. By the remarks above this is
enough. Suppose this is proved for all $n < m$ and we want to prove it for $m$.
Let $K$ be a finite extension of $F$ and $L$ a Galois extension of $K$ with
Galois group $G$ of order $m$. If $G$ has a proper normal subgroup $H$ then let
$K_{1}$ be $Fix(H)$. By inductive hypothesis, both $H^{2}(H,L^{*})$ and
$H^{2}(G/H,K_{1}^{*})$ are trivial. By Fact 4.3, so is $H^{2}(G,K^{*})$.
So we may assume $G$ to be simple. If $G$ is nonabelian then the order of $G$
is not a prime power. For each prime $p$ dividing the order of $G$, let $G_{p}$
be a $p$-Sylow subgroup of $G$ and let $K_{p}$ be the fixed field of $G_{p}$.
By induction hypothesis $H^{2}(G_{p},L^{*})$ is trivial, so by Fact 4.4 so is
$H^{2}(G,L^{*})$. So we may assume that $G$ is elementary abelian, of
cardinality $p$ say. Let $K_{1}$ be obtained from $K$ by adjoining all $pth$
roots of unity. Let $L_{1}$ be the compositum of $L$ and $K_{1}$. By Fact 4.3 it
is enough to prove that $H^{2}(Gal(L_{1}/K)$ is trivial. As all prime divisors
of the order of $Gal(K_{1}/K)$ are strictly less than $p$ it follows by the
induction hypothesis, as above,that $H^{2}(Gal(K_{1}/K),K_{1}^{*})$ is trivial.
Thus by Fact 4.3, it remains to prove only that $H^{2}(Gal(L_{1}/K_{1}),
L_{1}^{*})$ is trivial. Note that $L_{1}$ is an extension of $K_{1}$ of degree
$p$.
\newline
\noindent
So, changing notation, we are
reduced to showing that $H^{2}(Gal(L/K),L^{*})$ is trivial, when $G$ =
$Gal(L/K)$ has order $p$ and $K$, a finite extension of $F$, contains all $pth$
roots of unity. If $p$ is the characteristic,
then by Fact 4.2 (ii), the restriction of $N_{L/K}$ to $K^{*}$ is
already surjective, so $N_{L/K}:L^{*}\rightarrow K^{*}$ is
surjective, which by Fact 4.5 finishes the proof. So we may assume that $p$ is
prime to the characteristic. In that case, $L$ is a
Kummer extension of $K$ generated by a solution $\alpha$ to $x^{p} = a$ for
some $a\in K_{1}$. We will show that
$N_{L/K}:L^{*}\rightarrow K^{*}$ is surjective.

\noindent
Note that $\{1,\alpha, {\alpha}^{2},..,{\alpha}^{p-1}\}$ is a basis for $L$
over $K$, with respect to which any element of $L$ has coordinates
$x_{1},..,x_{p}$ from $K$. So the norm map can be represented as a map from
the set of $p$-tuples of elements of $K$ (not all zero) to $K^{*}$. Let $\omega$
be a primitive $pth$ root of unity. Then the conjugates of $\alpha$ under
$G$ are
$\alpha, \omega\alpha,..,{\omega}^{p-1}\alpha$. An easy computation shows that
for any $x_{1},x_{2}$ in $K$, $N_{L/K}(x_{1},x_{2},0,...0) = x_{1}^{p} +
ax_{2}^{p}$ if $p$ is not $2$ and $ = x_{1}^{p} - ax_{2}^{p}$ if $p=2$. Now $K$,
being a finite extension of
$F$ is also supersimple. By Fact 4.2, the multiplicative subgroup $T$ of $K^{*}$
of
$pth$ powers is definable and of finite index. By Proposition 3.4, both $T +
aT$ and $T - aT$ contain
$K^{*}$. Thus, the image of $N_{L/K}$ must be equal to $K^{*}$. By Fact 4.5,
$H^{2}(L/K) = 0$. The proof is complete.

\vspace{5mm}
\noindent
Note that Theorem 4.6 has the following consequence:
\begin{Corollary} Suppose that $F$ is a supersimple field, and $D$ is a
finite dimensional division algebra over $F$. Then $D$ is a field.
\end{Corollary}
{\em Proof.} By for example 15.8 of \cite{Lam}, if a division ring is
finite-dimensional over a subfield then it is finite-dimensional over its
centre.

\section{Supersimple division rings}
This section is devoted to a proof of:
\begin{Theorem} Any supersimple division ring is a field.
\end{Theorem}

Corollary 4.7 will play a crucial role in the proof.

We proceed to the proof of Theorem 5.1, which will go through various reductions
and cases.

We will assume that $D$ is a supersimple division ring (namely a division ring
type-definable in a big model $\bar{M}$ of a supersimple theory), which is not
commutative, and look for a contradiction. We make continuous use of Fact 4.1.
In particular $SU(D)$ is a monomial $\omega ^{\alpha}.n$ say ($\alpha$ an
ordinal and $n$ a positive integer).

Now for any finite subset $A$ of $D$, the centralizer of $A$ in $D$, $C_{D}(A)$
is a definable subdivision ring of $D$. Note that if $D_{1} < D_{2}$ are
infinite subdivision rings of $D$ then the index of $D_{1}$ in $D_{2}$
(additively) must be infinite, hence $SU(D_{1}) < SU(D_{2})$. It follows that
$D$ has the DCC on centralizers. Now the centre of $D$ is a field. So choosing a
smallest centralizer which is not a field, we may make:

\noindent
{\bf Assumption. For every $a\in D\setminus Z(D)$, $C_{D}(a)$ is a field.}

\vspace{5mm}
\begin{Lemma} $SU(Z(D)) < {\omega}^{\alpha}$ and also
$SU(C_{D}(a)) < {\omega}^{\alpha}$ for each $a\notin Z(D)$.
\end{Lemma}
{\em Proof.} Otherwise, $D$ has finite dimension over a subfield, contradicting
Corollary 4.7.

\vspace{5mm}
\begin{Lemma} (i) For any $a\notin Z(D)$ the conjugacy class of $a$ in $D$,
$a^{D}$ is generic, and moreover there are only finitely many such conjugacy
classes.
\newline
\noindent
(ii) Also for $a\notin Z(D)$, $\{x^{a}-x:x\in D\}$ is an additive subgroup of
finite index of $D$.
\end{Lemma}
{\em Proof.} (i) For $a\notin Z(D)$, $a^{D}$ is in definable bijection with
$D^{*}/C_{D}(a)^{*}$ which by Fact 4.1 and Lemma 5.2  has
$SU$-rank ${\omega}^{\alpha}.n$. Thus $a^{D}$ is generic in $D$.As the relation
of being in the same conjugacy class is an equivalence relation, there can be
at most finitely  many such generic classes.
\newline
\noindent
(ii) Note that $C_{D}(a)$ is precisely the kernel of the additive endomorphism
$\mu$ of $D$ which takes $x$ to $x^{a}-x$. So by Lemma 5.2, $SU(Ker(\mu)) <
\omega ^{\alpha}$. Fact 4.1 implies that $Im(\mu)$ has $SU$-rank
$\omega ^{\alpha}.n$, so has finite index.

\vspace{3mm}
\begin{Lemma} If $a,b \notin Z(D)$ do not commute, then $C_{D}(a)\cap C_{D}(b)
= Z(D)$.
\end{Lemma}
{\em Proof.} Let $c\in C_{D}(a)\cap C_{D}(b)$. If $c \notin Z(D)$ then by
the Assumption above, $C_{D}(c)$ is a field, so $a$ and $b$ commute,
contradiction.

\begin{Lemma} Assume that $char(D) = 0$. Then
\newline
(i) For any $b\in D\setminus Z(D)$ and $n>0$, $b^{n} \notin Z(D)$
\newline
\noindent
(ii) for all $a \in D\setminus Z(D)$, $a^{D}-a \supseteq Z(D)$.
\end{Lemma}
{\em Proof.} (i) Note that as $char(D) = 0$, $D$ has no additive subgroups of
finite index, and thus by Lemma 5.3 (ii), $\{x^{b}-x:x\in D\} = D$. Thus there
are nonzero $c\in Z(D)$ and $a\notin Z(D)$ such that $a^{b}-a = c$. Then $b$ is
in the normaliser in $D$ of $C_{D}(a)$ and $b \notin C_{D}(a)$. If $b^{n}\in
Z(D)$ then  the division ring generated by $C_{D}(a)$ and $b$ is noncommutative
and finite-dimensional over the field $C_{D}(a)$, contradicting Corollary 4.7
\newline
\noindent
(ii) As $a\notin Z(D)$, $(a+Z(D))\cap Z(D) = \emptyset$. Lemma 5.3 (i) says
there
are only finitely many noncentral conjugacy classes in $D$ and so one of them,
say $b^{D}$, must intersect $a+Z(D)$ in a generic subset of $a+Z(D)$ (namely a
subset of $a+Z(D)$ of maximal $SU$-rank). Without loss of generality, $b \in
a+Z(D)$. Thus $(b^{D}\cap (a+Z(D)) - b)$ which equals $(b^{D}-b)\cap Z(D)$ is a
generic subset of $Z(D)$ which is easily seen to be an additive subgroup, and
thus an additive subgroup of finite index. It follows, as we are in
characteristic $0$ that $b^{D}-b$ contains $Z(D)$. In particular, note $a\in
b+Z(D) \subset b^{D}$, so also $a^{D}-a$ contains $Z(D)$.

\vspace{5mm}
\noindent
With the above lemmas we can finish the proof. We divide into cases.

\vspace{5mm}
\noindent
{\bf Case I. $char(D) = p > 0$.}
\newline
\noindent
By Lemma 5.3(ii) and compactness there is a uniform bound, say $m$, on the index
in $D^{+}$ of the subgroup $\{x^{a}-x:x\in D\}$ for $a\in D\setminus Z(D)$. We
will show that every element of $D^{*}/Z(D)^{*}$ has finite order, which by
Kaplansky's Theorem (see 5.15 in \cite{Lam}) implies that $D$ is commutative. If
$C_{D}(a)$ has cardinality at most $m$ then clearly $a$ has exponent at most $m$
in $D$. Otherwise we can (by choice of $m$)
find nonzero $a' \in C_{D}(a)$ of the form $b^{a} - b$ for some $b\in D$. Note
that $a$ and $b$ do not commute. But $a^{p}$ and $b$ do commute
($b^{a} = b + a'$ implies $b^{a^{p}} = b + pa' = b$). So $a^{p}$ is in
$C_{D}(a)\cap C_{D}(b)$ which by Lemma 5.4 equals $Z(D)$. So $a$ has order $p$
in $D^{*}/Z(D)^{*}$.
This finishes the proof of Theorem 5.1 in the positive characteristic case.

\vspace{5mm}
\noindent
{\bf Case II. $char(D) = 0$.}
\newline
\noindent
Suppose to begin with that every element of $Z(D)$ has a square root in $Z(D)$.
Thus, we obtain infinitely many roots of unity $\{\omega _{i}:i < \omega\}$ in
$Z(D)$. Let $a\notin Z(D)$. Then for some $i<j$, $\omega _{i}.a$ and $\omega
_{j}.a$ are in the same conjugacy class (as by Lemma 5.3 (i) there are only
finitely many noncentral conjugacy classes). Thus $a$ is conjugate to $\omega
.a$ for some nontrivial root of unity $\omega$, say $a^{b} = \omega .a$. Note
then that $b \notin C_{D}(a)$, but $b^{m}\in C_{D}(a)$ for $m$ such that $\omega
^{m} = 1$. By Lemma 5.4, $b^{m} \in Z(D)$. This contradicts Lemma 5.5 (i).

\vspace{3mm}
\noindent
So we may assume that some $c\in Z(D)$ has no square root in $Z(D)$, and thus
by Lemma 5.5 (i) $c$ is not a square in  $D$. Choose $a \in D\setminus
Z(D)$, and
let $F$ = $C_{D}(a)$, a (supersimple) field containing $Z(D)$. Let $T$ =
$(F^{*})^{2}$, which is by supersimplicity a definable subgroup of $F^{*}$ of
finite index. Let $z \in Z(D)$ be nonzero. By Proposition 3.4 there are $b\in T$
and $e \in cT$, both generic over $z$ in $F$ such that $b - e = z$. Note that
$e$ is not a square in $D$ (as $c$ is not). As $e$ is generic in $F$, $e \notin
Z(D)$ so by Lemma 5.5 (ii), $e^{d} - e = z$ for some $d\in D$. Thus $e^{d} = b$.
But $b$ is a square, while $e$ and thus $e^{d}$ is not a square, contradiction.
This completes the proof of Theorem 5.1, and also the paper.

\end{document}